\documentclass[letterpaper]{article} 
\usepackage{aaai24}  
\usepackage{times}  
\usepackage{helvet}  
\usepackage{courier}  
\usepackage[hyphens]{url}  
\usepackage{graphicx} 
\usepackage{natbib}  
\usepackage{caption} 
\usepackage{algorithm}
\usepackage{algorithmic}

\usepackage{array}
\usepackage{subfig}
\usepackage{overpic}
\usepackage{anyfontsize}
\usepackage{graphicx}
\usepackage{multirow}
\usepackage{colortbl}
\usepackage{xcolor}
\usepackage{amsmath}
\usepackage{amssymb}
\usepackage{amsthm}
\usepackage{booktabs}
\usepackage{algorithm}
\usepackage{algorithmic}
\usepackage{threeparttable}
\usepackage[switch]{lineno}

\usepackage{newfloat}
\usepackage{listings}
\DeclareCaptionStyle{ruled}{labelfont=normalfont,labelsep=colon,strut=off} 
\floatstyle{ruled}
\newfloat{listing}{tb}{lst}{}
\floatname{listing}{Listing}
\title{A Reinforcement-Learning-Based \protect\\ Multiple-Column Selection Strategy for Column Generation}
\author{
    Haofeng Yuan, Lichang Fang, Shiji Song\thanks{Corresponding author.}\\
}
\affiliations{
    Department of Automation, BNRist, Tsinghua University\\
    \{yhf22, fanglc22\}@mails.tsinghua.edu.cn, shijis@mail.tsinghua.edu.cn
}

\begin{document}

\maketitle

\begin{abstract}
Column generation (CG) is one of the most successful approaches for solving large-scale linear programming (LP) problems. Given an LP with a prohibitively large number of variables (i.e., columns), the idea of CG is to explicitly consider only a subset of columns and iteratively add potential columns to improve the objective value. While adding the column with the most negative reduced cost can guarantee the convergence of CG, it has been shown that adding multiple columns per iteration rather than a single column can lead to faster convergence. However, it remains a challenge to design a multiple-column selection strategy to select the most promising columns from a large number of candidate columns. In this paper, we propose a novel reinforcement-learning-based (RL) multiple-column selection strategy. To the best of our knowledge, it is the first RL-based multiple-column selection strategy for CG. The effectiveness of our approach is evaluated on two sets of problems: the cutting stock problem and the graph coloring problem. Compared to several widely used single-column and multiple-column selection strategies, our RL-based multiple-column selection strategy leads to faster convergence and achieves remarkable reductions in the number of CG iterations and runtime.
\end{abstract}

\section{Introduction}

Column generation (CG) is a widely used approach for solving the linear programming (LP) relaxations of large-scale optimization problems that have a prohibitively large number of variables to deal with. It exploits the fact that the majority of feasible variables (i.e., columns) will not be part of an optimal solution. Therefore, CG starts with a subset of columns and gradually adds new columns that have the potential to improve the current solution, e.g., columns with a negative reduced cost (assuming a minimization problem), until no such columns exist and the current solution is proven optimal \cite{r_2}. CG is often combined with the branch-and-bound method to solve large-scale integer programming problems, which is called branch-and-price \cite{r_26}.

\begin{figure}[t]
\centering
\includegraphics[width=2.8in, keepaspectratio]{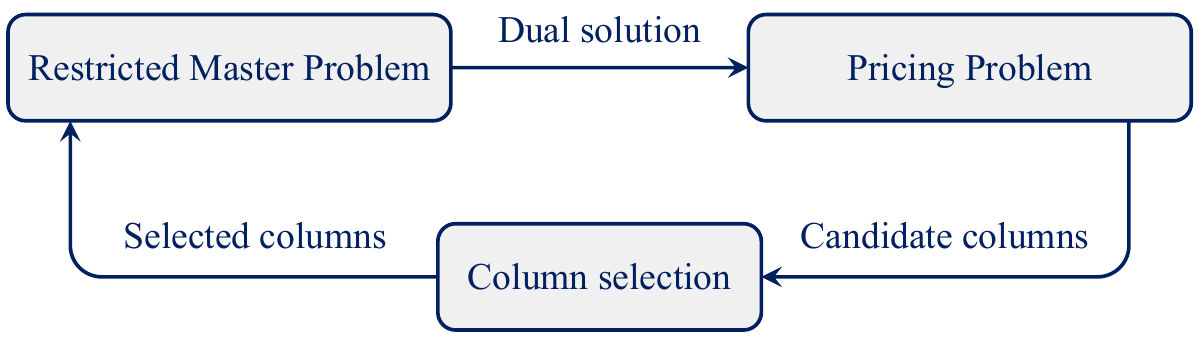}\\
\caption{The iterative process of CG.}
\label{fig_1}
\end{figure}

Specifically, CG follows an iterative process, as shown in Figure \ref{fig_1}. The original large-scale problem is decomposed into the restricted master problem (RMP) and the pricing problem (PP). CG starts by solving the RMP with a small subset of columns from the original problem. At each iteration, the RMP is solved using an LP solver (e.g., the simplex algorithm), and the dual solution is used to formulate the PP. The PP is a ``column generator'' that generates new columns (typically with negative reduced costs) to improve the current RMP solution. If such columns are found, they are added to the RMP to start a new iteration. Otherwise, it certifies the optimality of the current RMP solution for the original problem, and CG terminates. 

Typically, the column with the most negative reduced cost is selected to add to RMP at each iteration, which guarantees the convergence of CG and the optimality of the final solution \cite{r_1}. However, it often suffers from slow convergence, which limits its efficiency and usability. Previous research has shown that selecting multiple columns per iteration, including sub-optimal solutions of PP (even columns with non-negative reduced costs), can lead to faster convergence \cite{r_3}. This allows the RMP approximation to be improved, the optimal basis to be characterized faster, and hence to reduce the number of iterations. However, selecting multiple columns per iteration may result in a large fraction of useless columns that do not belong to the final optimal basis and increase the computation cost of RMP.

In general, the PP can generate a pool of feasible candidate columns, which are sorted according to the reduced cost, and the top-$k$ of them are greedily selected. In order to improve the selection, \citet{r_6} suggested that the RMP description can be improved by selecting non-correlated columns. Several diversification-based multiple-column selection strategies have been developed and shown to be effective in practice \cite{r_4,r_3}. However, despite the practical effectiveness of the diversification-based selection strategies, there is still no perfect column selection strategy proven to outperform or dominate the others.

Recently, reinforcement learning (RL) has shown impressive success in optimization tasks \cite{r_7,r_37,r_38}, which removes the need for substantial expert knowledge and pre-solved instances. In this paper, we propose a novel RL-based multiple-column selection strategy for CG. Specifically, we treat the iterative column selection in CG as a sequential decision task, and introduce an actor-critic style neural network that takes into account the column-constraint structure of RMP, the interrelations of candidate columns, and global properties of the problem instance. We use proximal policy optimization (PPO) \cite{r_8} to train the strategy to minimize the total number of iterations. Our RL-based multiple-column selection strategy is evaluated on two sets of problems: the cutting stock problem (CSP) \cite{r_9} and the graph coloring problem (GCP) \cite{r_10}. Experimental results demonstrate that our RL-based strategy outperforms several widely used single-column and multiple-column selection strategies in terms of the number of iterations and runtime. The main contributions of this paper can be summarized as follows:
\begin{itemize}
    \item We exploit RL to learn an effective multiple-column selection strategy for CG. To the best of our knowledge, it is the ﬁrst RL-based multiple-column selection strategy.
    \item We design an actor-critic style neural network that considers the column-constraint structure of RMP, the interrelations of candidate columns, and global properties of the problem instance, which allows to learn a column-relation-aware multiple-column selection strategy.
    \item We apply our approach to CSP and GCP, and experimental results show that it outperforms all baseline column selection strategies on various sizes of problems. Moreover, our RL-based framework can be easily applied to other problems solved based on CG.
\end{itemize}

\section{Related Work}

In this section, we review the acceleration methods for CG in the literature, with a focus on recent advances in machine learning (ML) techniques for column selection.

\paragraph{Acceleration Methods for Column Generation.}

Various techniques have been proposed in the literature to accelerate CG \cite{r_5}. One approach is to select ``better'' columns to add to RMP at each iteration. A classic approach is to add multiple columns rather than a single column with the most negative reduced cost. \citet{r_6} showed that the performance of CG is mathematically related to the variance-covariance matrix of selected columns: the convergence is accelerated with the selection of non-correlated columns. \citet{r_3} proposed two practical multiple-column selection strategies, which enhance the diversification of selected columns. Nevertheless, the effect of column selection is still not fully understood theoretically, and there is still no perfect column selection strategy proven to achieve the minimum number of iterations for CG. For faster convergence than existing hand-crafted column selection strategies, we apply RL to learn a multiple-column selection strategy that aims at minimizing the number of iterations through the interaction with the CG solution process.

Another approach is dual stabilization, which aims to form a ``better'' PP. For example, \citet{r_11} introduced a penalty function to reduce the oscillation of dual values. For a discussion on stabilization-based acceleration methods, please see \cite{r_12} and the references therein. More recently, \citet{r_35} proposed a learning-based method for predicting the optimal stabilization center of dual values in vehicle routing problems. We remark that our column selection strategy does not conflict with dual stabilization techniques and can be used synergistically for further improvement.

\paragraph{Machine-Learning-based Column Selection Strategy.}

Over the last few years, researchers have become increasingly interested in ML to accelerate optimization tasks \cite{r_28}, and several learning-based methods have been proposed for specific problems solved by CG \cite{r_15, r_14, r_13}. The closest works to ours are \cite{r_16} and \cite{r_17}, both leveraging ML for a better column selection strategy.

\citet{r_16} proposed a one-step lookahead ``expert'' to identify the columns that maximize the improvement of RMP in the next iteration, which is achieved by solving an extremely time-consuming mixed-integer linear programming (MILP). Then, they trained an ML model to cheaply imitate the decisions of the expensive MILP expert. They formulated the column selection procedure at each iteration as a classification task and trained the ML model in a supervised manner. The drawback of their approach is that the one-step lookahead expert is short-sighted because it only focuses on the very next iteration but disregards the interdependencies across iterations. Besides, it requires expensive pre-solved instances from previous executions of the MILP expert as training data, which may be unaffordable for large-scale applications. Moreover, the ML model only imitates the decision from the MILP expert, and thus it can never surpass the decisions for demonstration.

\citet{r_17} proposed a DQN-based single-column selection strategy that applies Q-learning to identify the ``best'' column at each iteration. They utilize the graph neural network (GNN) as a Q-function approximator to maximize the total expected future reward. While showing improved performance compared to the greedy single-column selection strategy, their framework can hardly be extended to multiple-column selection due to the exponential growth of action space and complex interdependencies in the column combinations, which limits its practical application (we implement a multiple-column variant of their DQN-based approach in our experiments). In contrast, our proposed neural network and learning scheme overcome these challenges and can derive an effective multiple-column selection strategy. Experiment results show that our RL-based multiple-column selection strategy outperforms the MILP expert used in \cite{r_16} and the DQN-based strategy proposed in \cite{r_17}.

\section{Basis of Column Generation} \label{sec_cg}

In this section, we use CSP as an example to introduce the mathematical formulation of CG. The CSP aims to determine the smallest number of rolls of length $L$ that have to be cut to satisfy the demands of $m$ customers, where customer $i$ demands $d_i$ pieces of orders of length $\ell_i, i=1,2,...,m$. \citet{r_9} proposed the CG formulation, in which the set $\mathcal{P}$ of all feasible cutting patterns is:
\begin{align*}
    \mathcal{P}=\left\{\left(\begin{array}{c}
    a_{1} \\
    \vdots \\
    a_{m}
    \end{array}\right) \in \mathbb{N}^{m} \mid \sum_{i=1}^{m} \ell_{i} a_{i} \leq L\right\}.
\end{align*}
Each pattern $p \in \mathcal{P}$ is denoted by a vector $\left( a_{1p}, \dots, a_{mp} \right)^{\mathrm{T}}$ $\in \mathbb{N}^{m}$, where $a_{ip}$ represents the number of pieces of length $\ell_i$ obtained in cutting pattern $p$. Let $\lambda_p$ be a decision variable that denotes the number of rolls cut using pattern $p \in \mathcal{P}$. The CSP is formulated as follows:
\begin{align*}
    \text{ min } & \sum_{p \in \mathcal{P}} \lambda_{p} \\
    \text { s.t. } & \sum_{p \in \mathcal{P}} a_{i p} \lambda_{p} \ge  d_{i}, \; i \in \left\{ 1,2,...,m \right\}, \\
    & \lambda_{p} \in \mathbb{N}, \; p \in \mathcal{P}.
\end{align*}
The objective function minimizes the total number of patterns used, equivalent to minimizing the number of rolls used. The $m$ constraints ensure all demands are satisfied.

This formulation usually has an extremely large number of decision variables as $\mathcal{P}$ is exponentially large. Therefore, the RMP is proposed for the linear relaxation with an initial set $\tilde{\mathcal{P}} \subset \mathcal{P}$. The RMP is defined as follows:
\begin{align*}
    \text{ min } & \sum_{p \in \tilde{\mathcal{P}}} \lambda_{p} \\
    \text { s.t. } & \sum_{p \in \tilde{\mathcal{P}}} a_{i p} \lambda_{p} \ge d_{i}, \; i \in \left\{ 1,2,...,m \right\}, \\
    & \lambda_{p} \ge 0, \; p \in \tilde{\mathcal{P}}.
\end{align*}
Let $u=\left(u_{1}, \ldots, u_{m}\right)^{\mathrm{T}}$ be the dual solution of the RMP. The columns that can potentially improve the solution of RMP are given by the solution to the following knapsack problem, which is referred to as the PP:
\begin{align*}
    \text{ max } & \sum_{i=1}^{m} u_{i}a_i \\
    \text { s.t. } & \sum_{i=1}^{m} \ell_{i} a_{i} \leq L, \; i \in \left\{ 1,2,...,m \right\}, \\
    & a_{i} \in \mathbb{N}, \; i \in \left\{ 1,2,...,m \right\}.
\end{align*}
The PP generates feasible patterns (columns), represented as vector $\left(a_{1p}, \dots, a_{mp} \right)^{\mathrm{T}}$, to be added to $\tilde{\mathcal{P}}$ for the next iteration. In general, several sub-optimal solutions of PP, which form a candidate column pool, can be obtained through dynamic programming methods or commercial solvers such as \emph{Gurobi} \cite{r_36}. We can select one or multiple columns from the candidate column pool to add to RMP for the next iteration (see Figure \ref{fig_1}).

\section{Methodology} \label{sec_method}

In this section, we present the details of our RL-based multiple-column selection strategy.

\subsection{MDP Formulations}

We treat the CG solution process as the environment for the RL agent. We formulate the multiple-column selection task for CG as a Markov decision process (MDP):

\paragraph{State $\mathcal{S}$.} \label{sec_state}

\begin{figure}[b]
\centering
\includegraphics[width=2.6in, keepaspectratio]{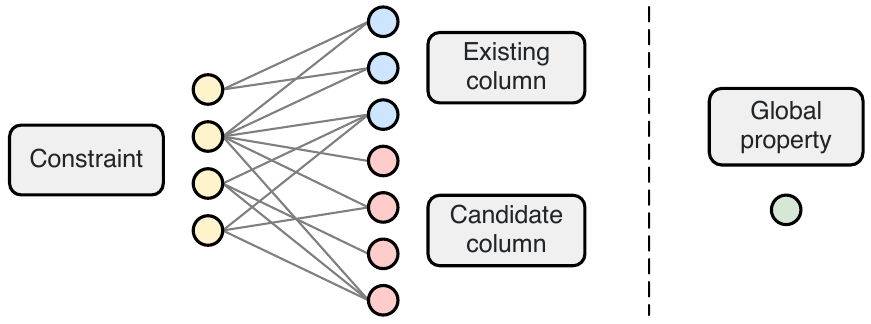}\\
\caption{A toy example of state. The left part illustrates the bipartite graph representation of the current RMP, including 4 constraint nodes, 3 existing column nodes, and 4 candidate column nodes. The right part denotes the global feature vector for the problem instance.}
\label{fig_2}
\end{figure}

The state describes the information about the current CG status. which is provided for the RL agent. As illustrated in Figure \ref{fig_2}, the state is defined to include 1) a bipartite graph representation of the current RMP and candidate columns, and 2) global properties of the problem instance.

As introduced in \cite{r_18}, an LP can be represented as a bipartite graph with constraint nodes $\mathcal{C}$ and column nodes $\mathcal{V}$. We further incorporate candidate columns into the bipartite graph representation, where column nodes are divided into existing columns in the current RMP and candidate columns to be selected (blue nodes and red nodes in Figure \ref{fig_2}). An edge $\left( v, c \right)$ exists between a node $v \in \mathcal{V}$ and a node $c \in \mathcal{C}$ if column $v$ contributes to constraint $c$. State information corresponding to the columns (e.g., solution value, reduced cost) and constraints (e.g., slack, dual value) are represented as node features. In addition to the bipartite graph representation, we represent the properties associated with the problem instance (e.g., the number of constraints, maximum constraint coefﬁcient) as an additional global feature vector (the green node in Figure \ref{fig_2}).

\paragraph{Action $\mathcal{A}$.} \label{sec_action}

At each iteration, we select $k$ columns from the pool of $n$ candidate columns generated by the PP, and add them to RMP for the next iteration. The action space $\mathcal{A}$ contains all possible $k$-combinations of the $n$ candidate columns, i.e., $\left | \mathcal{A} \right | = \binom{n}{k}$. The RL agent returns a probability distribution over the action space, and we sample an action from that distribution, which can be seen as a multiple-column selection strategy.

\paragraph{Transition $\mathcal{T}$.}

The transition rule is deterministic. Once an action is selected, the corresponding $k$ columns are added to the RMP to start a new iteration.

\paragraph{Reward $\mathcal{R}$.}

The goal of the RL agent is to minimize the total number of iterations. We design a reward function consisting of a unit penalty for each additional iteration and two auxiliary components: 1) the decrease in the objective value of RMP, and 2) the sum of cosine distances of selected columns, which is inspired by the observation in \cite{r_4}. The immediate reward at time step $t$ is:
\begin{align*}
\text{\fontsize{7.8pt}{7.8pt}\selectfont
$r_t = -1 + \alpha \cdot \left(\frac{\text{obj}_{t-1}-\text{obj}_{t}}{\text{obj}_{0}}\right) + \beta \cdot \sum_{u_i, u_j \in \mathcal{C}_s} \left( 1 - \frac{\left \langle u_i, u_j \right \rangle}{\left \| u_i \right \| \cdot  \left \| u_j \right \|} \right),$}
\end{align*}
\normalsize
where ($\text{obj}_{t-1}-\text{obj}_{t}$) is the decrease in the objective value, normalized by the objective value $\text{obj}_{0}$ of the initial RMP; $\mathcal{C}_s$ is the set of selected columns, and $\left( 1 - \frac{\left \langle u_i, u_j \right \rangle}{\left \| u_i \right \| \cdot  \left \| u_j \right \|} \right)$ is the cosine distance between column vectors $u_i$ and $u_j$; $\alpha$ and $\beta$ are non-negative weight hyperparameters.

\subsection{Model} \label{sec_model}

We use PPO \cite{r_8} as the training algorithm for our multiple-column selection strategy. PPO is a deep reinforcement learning algorithm based on the actor-critic architecture. Given a state $s$, the critic estimates the value function $V(s)$, and the actor gives a probability distribution $\pi = \left(\pi\left(a_1 \mid s\right), \pi\left(a_2\mid s\right), \dots, \pi\left(a_{\left | \mathcal{A} \right |}\mid s\right)\right)$ over the action space $\mathcal{A}$. An action is sampled from the probability distribution, and the corresponding $k$ columns are selected and added to the RMP to start a new iteration.

We propose an actor-critic style neural network for the RL-based multiple-column selection strategy. The network consists of three components: an encoder, a critic decoder, and an actor decoder. The details of the neural network architecture are described below:

\paragraph{Encoder.}

\begin{figure}[t]
\centering
\includegraphics[width=2.5in, keepaspectratio]{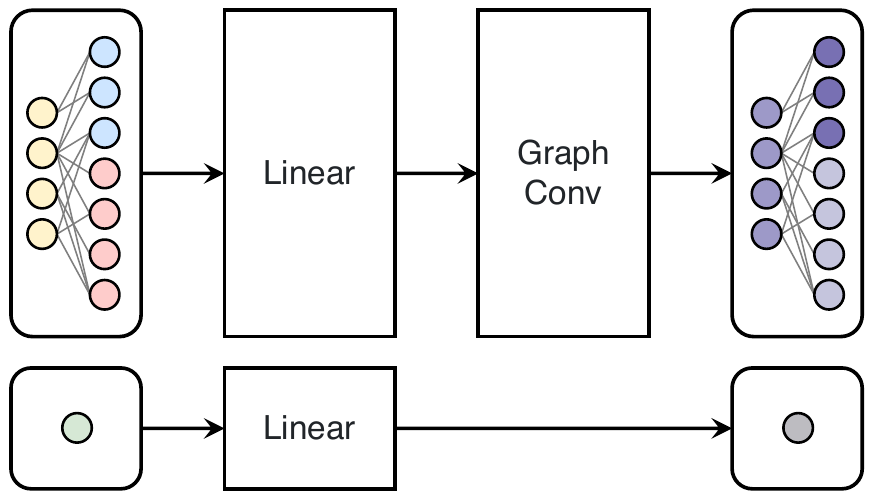}\\
\caption{The architecture of the encoder. Colored nodes denote feature vectors or embeddings.}
\label{fig_3}
\end{figure}

As introduced above, the state is represented by a bipartite graph and a global feature vector. The encoder takes the bipartite graph and the global feature vector as input to produce the embeddings for the current state. The architecture of the encoder is shown in Figure \ref{fig_3}.

Specifically, the bipartite graph and the global feature vector are embedded separately.
For the bipartite graph, the encoder first computes the initial node embeddings from raw node features through a learned linear projection. Then, the node embeddings of the bipartite graph are updated through $N_1$ graph convolutional layers. Each layer proceeds in two phases: the first phase is performed to update the constraint node embeddings, followed by the second phase that updates the column node embeddings. Both phases are implemented using the graph isomorphism network (GIN) \cite{r_19} with residual connections. Let $x_{c}^{(\ell)}$ and $x_{v}^{(\ell)}$ denote the embeddings for constraint node $c \in \mathcal{C}$ and column node $v \in \mathcal{V}$ at layer $\ell$. The node embeddings are updated as follows:
\begin{align*}
    & \text{\fontsize{9pt}{9pt}\selectfont $x_{c}^{(\ell)}=\operatorname{MLP}^{(\ell)}_{\mathcal{C}}\left(\left(1+\epsilon \right) \cdot x_{c}^{(\ell-1)}+\sum_{v_i \in \mathcal{N}\left(c\right)} x_{v_i}^{(\ell-1)}\right) + x_{c}^{(\ell-1)},$} \\
    & \text{\fontsize{9pt}{9pt}\selectfont $x_{v}^{(\ell)}=\operatorname{MLP}^{(\ell)}_{\mathcal{V}}\left(\left(1+\epsilon \right) \cdot x_{v}^{(\ell-1)}+\sum_{c_i \in \mathcal{N}\left(v\right)} x_{c_i}^{(\ell)}\right) + x_{v}^{(\ell-1)},$}
\end{align*}
where $\operatorname{MLP}^{(\ell)}_{\mathcal{C}}$ and $\operatorname{MLP}^{(\ell)}_{\mathcal{V}}$ are multi-layer perceptrons (MLPs) for updating the constraint node embeddings and column node embeddings, respectively, and $\mathcal{N}\left(v\right)$ denotes the neighborhood set of node $v$.

\begin{figure*}[t]
\centering
\subfloat[The critic decoder]{\includegraphics[width=2.625in, keepaspectratio]{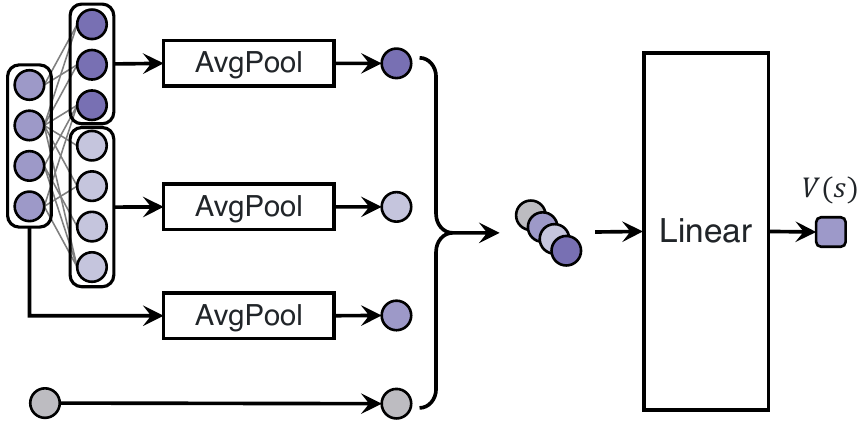}}\hspace{5pt}
\subfloat[The actor decoder]{\includegraphics[width=4.2in, keepaspectratio]{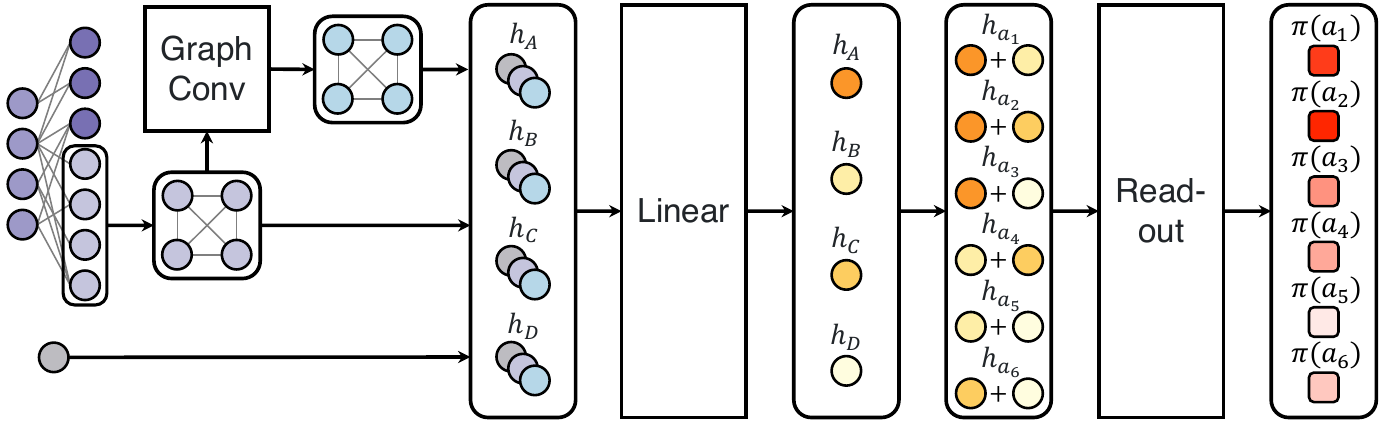}}
\caption{The architecture of the critic decoder and actor decoder. Here shows a toy example of selecting 2 from 4 candidates.}
\label{fig_5}
\end{figure*}

The global feature vector is embedded through $N_2$ linear layers, each followed by a LeakyReLU activation function.

\paragraph{Critic Decoder.}

The critic decoder maps the latent embeddings of state $s$ into the estimated value function $V(s)$. The architecture of the critic decoder is illustrated in Figure \ref{fig_5}(a). In the critic decoder, the node embedding vectors associated with the existing columns, candidate columns, and constraints are respectively pooled by average, and then concatenated together with the embedding vector of global features, which contains information of the current RMP as well as the global properties of the problem instance. Then, the concatenated vector is passed through an $N_3$ layer MLP to estimate the value function $V(s)$.

\paragraph{Actor Decoder.}

Based on the embeddings produced by the encoder, the action decoder outputs a probability distribution $\pi = \left(\pi\left(a_1 \mid s\right), \pi\left(a_2\mid s\right), \dots, \pi\left(a_{\left | \mathcal{A} \right |}\mid s\right)\right)$ over the action space. An action is sampled from the probability distribution, determining which $k$ columns are selected by the RL-based multiple-column selection strategy.

The interrelations, especially the similarity between candidate columns, are crucial to the multiple-column selection. Therefore, we explicitly model the message-passing between candidate columns. We first create a complete graph, with each node corresponding to a candidate column. The node embeddings of the complete graph are initialized as the final embeddings of candidate column nodes from the bipartite graph. We associate the distance (e.g., Jaccard distance, cosine distance) between candidate column vectors as the initial edge features. As shown in Figure \ref{fig_5}(b), the candidate column node embeddings of the bipartite graph are used to create the complete graph.

We apply the graph attention network (GAT) with edge features \cite{r_20,r_21} to update the embeddings of the complete graph through message-passing between candidate columns. Then, for each candidate column, we concatenate its node embedding from the bipartite graph, its node embedding from the complete graph, and the global embedding together. The concatenated embedding for each candidate column is processed through an $N_4$ layer MLP to obtain the final embedding vector ($h_A$ to $h_D$ in Figure \ref{fig_5}(b)).

For each action $a_i$, we define its representation vector $h_{a_i}$ $= \sum_{v_j \in \mathcal{C}_s\left(a_i\right)}h_{v_j}$, where $\mathcal{C}_s\left(a_i\right)$ denotes the set of columns selected in action $a_i$, and $h_{v_j}$ is the final embedding vector of candidate column $v_j$. The probability to select $a_i$ is computed through a learnable nonlinear readout function:
\begin{align*}
    \pi\left(a_i \mid s\right) = \operatorname{softmax} \left(C \cdot \operatorname{tanh}\left(w_o^{\mathrm{T}} \cdot \operatorname{ReLU}\left(W_o \cdot h_{a_i}\right)\right)\right),
\end{align*}
where $w_o$ and $W_o$ are learnable vector and matrix, respectively, and $C$ is the clipping coefficient ($C=10$). The actor acts as a multiple-column selection strategy, which takes the current state as input and outputs a probability distribution over the action space.

Note that we are using a learnable nonlinear mapping to derive the probability distribution for each $k$-combination of columns. This is by no means a simple addition of individual scores for the corresponding $k$ columns. An action is sampled from the probability distribution, and the $k$ candidate columns in the corresponding combination are added to the RMP for the next iteration.

\section{Evaluation}

We evaluate our proposed RL-based multiple-column selection strategy on two sets of problems: the CSP and the GCP. Both problems are well-known for the linear relaxation effectively solved using CG. Experimental results demonstrate that our RL-based strategy outperforms several widely used single-column selection strategies and multiple-column selection strategies.

\subsection{Experiment Task}

\paragraph{Cutting Stock Problem.}

The CG formulation of CSP has been introduced in the previous section. The problem instances are generated according to the rules for random instances in BPPLIB \cite{r_24, r_23}, a widely used benchmark for bin packing and cutting stock problems. We divide the CSP instances into three categories: easy, normal, and hard, corresponding to the roll length $L=50, 100,$ and $200$. We generated 1000, 200, and 100 instances of the three instance categories for evaluation, respectively. Instances for training are randomly generated and solved using CG as the environment for the RL agent.

\paragraph{Graph Coloring Problem.}
The GCP aims to assign a minimum number of colors to the nodes on a graph, such that every pair of adjacent nodes does not share the same color \cite{r_22}. In the CG formulation for GCP, the RMP can be expressed as using a minimum number of maximal independent sets (MISs) to cover all the nodes, while the PP is modeled as the maximum weight independent set problem (MWISP) to search for new MISs with the most negative reduced cost. The details are described in \cite{r_10}. Similar to the CSP, the GCP instances are divided into three categories, corresponding to the number of nodes $N=30, 40,$ and $50$ respectively. The CGP instances are generated according to the rules for random graphs in \cite{r_10}.

\subsection{Hyperparameter Conﬁguration}
We implement our model to learn the RL-based multiple-column strategy. We select 5 out of 10 candidate columns at each iteration, which strikes a balance between the number of iterations and the cost per iteration in our task. To guarantee convergence, we force the optimal solution of the PP to always be selected. We set the number of layers $N_1=N_2=N_3=N_4=3$ for the MLPs in the encoder and decoder. The weights in the reward function are set to $\alpha=300$ and $\beta=0.02$ to balance the reward scales and the discount factor $\gamma$ is set to $0.9$. We use PPO with a clipping threshold $\epsilon=0.2$, and the Adam optimizer with a learning rate $1 \times 10^{-3}$ to train the RL model. The hyperparameter configuration is fixed across all instance categories of CSP and GCP.

\subsection{Comparison Evaluation}
We compare our RL-based multiple-column strategy with several well-established single-column and multiple-column selection strategies, as well as the multiple-column selection strategy using the MILP expert proposed in \cite{r_16} and the DQN-based approach proposed in \cite{r_17}. The details of baseline strategies for comparison are as follows:

\paragraph{Single-column selection strategy:}
\begin{itemize}
    \item  \emph{Greedy single-column selection (Greedy-S)}: Always select the column with the most negative reduced cost.
    \item  \emph{Random single-column selection (Random-S)}: Randomly select a column from the candidate column pool.
    \item  \emph{DQN-based single-column selection (DQN-S)}: Selection strategy based on DQN in \cite{r_17}.
\end{itemize}
\paragraph{Multiple-column selection strategy:}
\begin{itemize}
    \item  \emph{Greedy multiple-column selection (Greedy-M)}: Always select the top-$k$ columns according to the reduced cost.
    \item  \emph{Random multiple-column selection (Random-M)}: Randomly select $k$ columns from the candidate column pool.
    \item  \emph{MILP expert (MILP-M)}: Selection strategy using the MILP expert in \cite{r_16}.
    \item  \emph{Diversification-based column selection (Diverse-M)}: A modified strategy of CGDS \cite{r_3} to fit our task: We first sort the candidate columns by their reduced costs, and prioritize the candidate columns that are disjoint from the already selected columns, if there exists one in the remaining pool.
\end{itemize}

\begin{table}[t]
\renewcommand\arraystretch{1.1}
    \setlength{\tabcolsep}{3.5pt}
    \centering
    \fontsize{9pt}{9pt}\selectfont
    \begin{tabular}{l|cc|cc|cc}
        \toprule
        \multirow{2.5}*{Strategy} & \multicolumn{2}{c}{CSP (Easy)} & \multicolumn{2}{c}{CSP (Normal)} & \multicolumn{2}{c}{CSP (Hard)}\\
        \cmidrule(lr){2-3} \cmidrule(lr){4-5} \cmidrule(lr){6-7}
         & \# Itr & Time & \# Itr & Time & \# Itr & Time \\
        \midrule
        \midrule
        Greedy-S  & 37.68 & 228.92 & 89.20 & 186.78 & 171.44 & 301.07\\
        Random-S  & 62.63 & 374.80 & 116.82 & 257.04 & 205.31 & 376.16\\
        DQN-S     & 35.54 & 215.65 & 88.95  & 178.52 & / & /  \\
        \midrule
        Greedy-M  & 12.03 & 75.34 & 27.01 & 62.60 & 52.13 & 96.42 \\
        Random-M  & 13.97 & 84.83 & 28.26 & 63.78 & 51.43 & 96.51 \\
        MILP-M    & 10.65 & 96.23 & 23.46 & 81.17 & 44.99 & 147.19 \\
        Diverse-M & 11.45 & 74.59 & 25.03 & 59.32 & 47.24 & 93.60 \\
        Ours & \textbf{10.33} & \textbf{67.84} & \textbf{22.85} & \textbf{55.05} & \textbf{43.95} & \textbf{87.95} \\
        \bottomrule
    \end{tabular}
    \caption{Comparison results on the CSP, in terms of the average number of iterations per instance and total runtime (in seconds) over the evaluation instance set.}
    \label{table_1}
\end{table}

For a fair comparison, we set the same number of candidate columns and the columns to select in all multiple-column selection strategies. The candidate columns are generated as the 10 columns with the most negative reduced cost from PP. We report the evaluation metrics of 1) the average number of iterations per instance and 2) the total runtime over the evaluation set. Generally, all heuristic strategies and the RL-based strategies (on GPU), except MILP-M, require negligible time for a selection decision, so the comparison of the average number of iterations is approximately equivalent to the comparison of runtime for these strategies. We remark that MILP-M is practically intractable because the time taken by the MILP expert in MILP-M is even larger than the time for RMP and PP \cite{r_16}. As shown in the experimental results, even if MILP-M requires fewer iterations for convergence compared to other heuristic column selection strategies, its runtime is significantly larger due to the expensive decisions from the MILP expert.

\begin{table}[t]
\renewcommand\arraystretch{1.1}
    \setlength{\tabcolsep}{4.5pt}
    \centering
    \fontsize{9pt}{9pt}\selectfont
    \begin{tabular}{l|cc|cc|cc}
        \toprule
        \multirow{2.5}*{Strategy} & \multicolumn{2}{c}{GCP (Easy)} & \multicolumn{2}{c}{GCP (Normal)} & \multicolumn{2}{c}{GCP (Hard)}\\
        \cmidrule(lr){2-3} \cmidrule(lr){4-5} \cmidrule(lr){6-7}
         & \# Itr & Time & \# Itr & Time & \# Itr & Time \\
        \midrule
        \midrule
        Greedy-M  & 18.30 & 75.16 & 29.00 & 99.54 & 39.04 & 123.23\\
        Random-M  & 19.04 & 78.81 & 30.91 & 104.79 & 41.01 & 125.26 \\
        MILP-M    & 17.07 & 93.78 & 25.86 & 128.71 & 34.19 & 188.44 \\
        Diverse-M & 18.07 & 74.36 & 28.17 & 96.14  & 37.78 & 119.90 \\
        Ours & \textbf{15.19} & \textbf{62.06} & \textbf{24.93} & \textbf{87.31} & \textbf{34.11} & \textbf{111.56}  \\
        \bottomrule
    \end{tabular}
    \caption{Comparison results on the GCP.}
    \label{table_2}
\end{table}

\paragraph{Results on CSP.}
The comparison results on CSP are reported in Table \ref{table_1}. All the multiple-column selection strategies achieve significantly faster convergence than the single-column selection strategies, especially on large-scale problem instances. Diverse-M requires the least runtime among the baseline strategies. While MILP-M requires fewer iterations than Diverse-M, it takes a much larger runtime due to the extremely expensive column selection decisions.

The experimental results show that our RL approach can learn a stronger strategy for column selection, which is implicitly represented by the neural network. The RL-based multiple-column selection strategy outperforms all baseline column selection strategies in the three instance categories of various sizes. Compared to the best baseline strategy (Diverse-M), our RL-based multiple-column selection strategy yields a total runtime reduction of \textbf{9.05\%}, \textbf{7.20\%}, and \textbf{6.04\%} in the three instance categories, respectively. In addition, it is worth mentioning that our RL-based multiple-column selection strategy requires even fewer iterations than MILP-M. This is mainly because the MILP expert focuses only on the current step and its goal is to minimize the objective value for the very next iteration, whereas our RL agent aims to minimize the total number of iterations and takes into consideration the long-term effect of currently selected columns on future iterations.

\paragraph{Results on GCP.}
Since the single-column selection strategies have been demonstrated to be ineffective on CSP, we conduct experiments on GCP using only multiple-column selection strategies. As reported in Table \ref{table_2}, it shows similar results to the experiments on CSP. Compared to Diverse-M, our RL-based strategy reduces the total runtime by \textbf{16.54\%}, \textbf{9.18\%}, and \textbf{6.96\%} in the three instance categories, respectively. The evaluation results on CSP and GCP demonstrate the effectiveness of our RL-based multiple-column selection strategy on different types of problems.

\begin{table*}[t]
\renewcommand\arraystretch{1.1}
    \setlength{\tabcolsep}{8pt}
    \centering
    \fontsize{9pt}{9pt}\selectfont
    \begin{tabular}{l|cc|cc|cc|cc}
        \toprule
        \multirow{2.5}*{Strategy} & \multicolumn{2}{c}{CSP ($L$=500)} & \multicolumn{2}{c}{CSP ($L$=1000)} & \multicolumn{2}{c}{GCP ($N$=75)} & \multicolumn{2}{c}{GCP ($N$=100)}\\
        \cmidrule(lr){2-3} \cmidrule(lr){4-5} \cmidrule(lr){6-7} \cmidrule(lr){8-9}
         & \# Itr & Time & \# Itr & Time & \# Itr & Time & \# Itr & Time \\
        \midrule
        \midrule
        Greedy-M  & 78.80 & 143.25 & 98.82 & 221.75 & 60.88 & 287.52 & 83.60 & 449.95 \\
        Random-M  & 74.96 & 136.66 & 85.74 & 206.81 & 65.50 & 301.77 & 87.73 & 469.44 \\
        MILP-M    & 67.26 & 247.00 & \textbf{81.92} & 482.13 & 53.40 & 438.38 & 75.22 & 788.65 \\
        Diverse-M & 70.44 & 130.84 & 83.74 & 196.92 & 57.80 & 259.46 & 83.07 & 449.74 \\
        Ours & \textbf{66.94} & \textbf{118.89} & 82.04 & \textbf{183.76} & \textbf{52.96} & \textbf{246.63} & \textbf{74.80} & \textbf{398.37} \\
        \bottomrule
    \end{tabular}
    \caption{Generalization performance of our RL-based multiple-column selection strategy to larger problem sizes.}
    \label{table_3}
\end{table*}

\begin{table}[t]
    \renewcommand\arraystretch{1.1}
    \setlength{\tabcolsep}{3.5pt}
    \centering
    \begin{threeparttable}
    \fontsize{9pt}{9pt}\selectfont
    \begin{tabular}{l|cc|cc|cc}
        \toprule
        \multirow{2.5}*{Strategy} & \multicolumn{2}{c}{CSP (Easy)} & \multicolumn{2}{c}{CSP (Normal)} & \multicolumn{2}{c}{CSP (Hard)}\\
        \cmidrule(lr){2-3} \cmidrule(lr){4-5} \cmidrule(lr){6-7}
         & \# Itr & Time & \# Itr & Time & \# Itr & Time \\
        \midrule
        \midrule
        Greedy-M  & 12.03 & 75.34 & 27.01 & 62.60 & 52.13 & 96.42 \\
        Random-M  & 13.97 & 84.83 & 28.26 & 63.78 & 51.43 & 96.51 \\
        Diverse-M & 11.45 & 74.59 & 25.03 & 59.32 & 47.24 & 93.60 \\
        \midrule
        DQN-M & 11.50 & 74.71 & 25.94 & 60.36 & 48.59 & 94.51 \\
        Variant 1\tnote{\romannumeral1} & 10.91 & 71.48 & 24.38 & 58.40 & 45.75 & 92.16 \\
        Variant 2\tnote{\romannumeral2} & 10.40 & 68.21 & 23.03 & 56.87 & 44.92 & 89.30 \\
        Complete Model & \textbf{10.33} & \textbf{67.84} & \textbf{22.85} & \textbf{55.05} & \textbf{43.95} & \textbf{87.95} \\
        \bottomrule
    \end{tabular}
    \begin{tablenotes}
        \vspace{0.3em}
        \item[\romannumeral1] with the embeddings of the complete graph removed.
        \item[\romannumeral2] with the embedding of the input global features removed.
    \end{tablenotes}
    \end{threeparttable}
    \caption{Performance of the RL-based multiple-column selection strategies using different neural network architectures.}
    \label{table_4}
\end{table}

\subsection{Generalization Evaluation}

Generalization across different instance sizes is a highly desirable property for learning-based models. The ability to generalize across sizes would allow the RL-based strategy to scale up to larger instances while training more efficiently on smaller instances. Table \ref{table_3} presents the generalization performance of our RL-based multiple-column selection strategy, which is trained on instances from the hard category but evaluated on much larger instances. For the CSP, the model is trained on instances with $L=200$ and evaluated on instances with $L=500$ and $L=1000$; for the GCP, the model is trained on instances with $N=50$ and evaluated on instances with $N=75$ and $N=100$. Our RL-based multiple-column selection strategy still shows advantages over most baseline column selection strategies on the evaluation instances, which are several times larger than the training instances. It is demonstrated that our RL-based multiple-column selection strategy has learned useful and effective selection principles that are invariant to the size of problem instances.

\subsection{Ablation Study} \label{Ablation Study}

We have demonstrated the effectiveness of our proposed RL-based multiple-column selection strategy on CSP and GCP. To show the effect of different components of the neural network, we consider two variants of the complete model: 1) removing the embeddings of the complete graph and 2) removing the embedding of input global features. Other components remain unchanged.

We conduct the ablation evaluation on CSP and the results are reported in Table \ref{table_4}. Compared to the complete model, the performance of both variants is degraded on all three instance categories. The results show that the embeddings of the complete graph and the embeddings of explicit global features both provide benefits for the learned multiple-column selection strategy. Notably, the performance of the learned multiple-column selection strategy decreases obviously when the embeddings of the complete graph are removed, which further highlights the importance of the candidate column interrelations in column selection.

To show the advantage of our approach over the framework of \cite{r_17}, we conduct an extended experiment using a modified implementation of their DQN-based approach for multiple-column selection (DQN-M in Table \ref{table_4}), where we select the top-$k$ columns based on their Q-values. DQN-M outperforms the greedy and random baselines but underperforms Diverse-M. This is because DQN-M independently selects $k$ columns with higher Q-values: while each of them can individually lead to good convergence, the combination of them is not necessarily the best, just as 5 Michael Jordans in a basketball team may not be better than a reasonable lineup. In other words, DQN-M is selecting ``the top-$k$ columns'', while our RL-based multiple-column selection strategy is devoted to selecting ``the best $k$-combination''.

\section{Conclusion}

In this paper, we propose an RL-based multiple-column selection strategy for CG. We formulate the multiple-column selection task as an MDP, and introduce an actor-critic style neural network that takes into account the column-constraint structure of RMP, the interrelations of candidate columns, as well as global properties of the problem instance. We evaluate our proposed RL-based multiple-column selection strategy on two sets of problems: the CSP and the GCP. Experimental results show that our RL-based multiple-column selection strategy outperforms all baseline single-column and multiple-column selection strategies in the three instance categories of various sizes. Extensive experiments also demonstrate the ability of our RL-based model to generalize to larger-scale problem instances.

Despite the significant performance of the RL-based multiple-column selection strategy, more progress can be made in exploring column selection strategies to select different numbers of columns adaptively based on the problem properties and solution stages. It is challenging for hand-crafted rules due to the various characteristics of different types and sizes of problems, but may be learned using an RL agent through the interaction with the CG solution environment. In addition, how to incorporate the learning-based column selection strategy with other acceleration methods for CG, such as dual stabilization and PP reduction, is also a possible direction for future efforts.

\section*{Acknowledgments}
This work was supported in part by the National Natural Science Foundation of China under Grant 61936009; and in part by the National Science and Technology Innovation 2030 Major Project of the Ministry of Science and Technology of China under Grant 2018AAA0101604. We also thank Wanlu Yang and Peng Jiang for their valuable comments and corrections.

\bibliography{reference}

\begin{thebibliography}{30}
\providecommand{\natexlab}[1]{#1}

\bibitem[{Babaki, Jena, and Charlin(2021)}]{r_35}
Babaki, B.; Jena, S.~D.; and Charlin, L. 2021.
\newblock Neural column generation for capacitated vehicle routing.
\newblock In \emph{AAAI-22 Workshop on Machine Learning for Operations
  Research}.

\bibitem[{Barnhart et~al.(1998)Barnhart, Johnson, Nemhauser, Savelsbergh, and
  Vance}]{r_26}
Barnhart, C.; Johnson, E.~L.; Nemhauser, G.~L.; Savelsbergh, M. W.~P.; and
  Vance, P.~H. 1998.
\newblock Branch-and-price: Column generation for solving huge integer
  programs.
\newblock \emph{Operations Research}, 46(3): 316--329.

\bibitem[{Bengio, Lodi, and Prouvost(2021)}]{r_28}
Bengio, Y.; Lodi, A.; and Prouvost, A. 2021.
\newblock Machine learning for combinatorial optimization: {A} methodological
  tour d'horizon.
\newblock \emph{European Journal of Operational Research}, 290(2): 405--421.

\bibitem[{Chi et~al.(2022)Chi, Aboussalah, Khalil, Wang, and
  Sherkat-Masoumi}]{r_17}
Chi, C.; Aboussalah, A.~M.; Khalil, E.~B.; Wang, J.; and Sherkat-Masoumi, Z.
  2022.
\newblock A deep reinforcement learning framework for column generation.
\newblock In \emph{Advances in Neural Information Processing Systems}.

\bibitem[{Delorme, Iori, and Martello(2016)}]{r_24}
Delorme, M.; Iori, M.; and Martello, S. 2016.
\newblock Bin packing and cutting stock problems: Mathematical models and exact
  algorithms.
\newblock \emph{European Journal of Operational Research}, 255(1): 1--20.

\bibitem[{Delorme, Iori, and Martello(2018)}]{r_23}
Delorme, M.; Iori, M.; and Martello, S. 2018.
\newblock {BPPLIB:} A library for bin packing and cutting stock problems.
\newblock \emph{Optimization Letters}, 12(2): 235--250.

\bibitem[{Desaulniers, Desrosiers, and Solomon(2002)}]{r_5}
Desaulniers, G.; Desrosiers, J.; and Solomon, M.~M. 2002.
\newblock Accelerating strategies in column generation methods for vehicle
  routing and crew scheduling problems.
\newblock \emph{Essays and Surveys in Metaheuristics}, 309--324.

\bibitem[{Desaulniers, Desrosiers, and Solomon(2006)}]{r_2}
Desaulniers, G.; Desrosiers, J.; and Solomon, M.~M. 2006.
\newblock \emph{Column generation}.
\newblock Springer Science \& Business Media.

\bibitem[{du~Merle et~al.(1999)du~Merle, Villeneuve, Desrosiers, and
  Hansen}]{r_11}
du~Merle, O.; Villeneuve, D.; Desrosiers, J.; and Hansen, P. 1999.
\newblock Stabilized column generation.
\newblock \emph{Discrete Mathematics}, 194(1-3): 229--237.

\bibitem[{Gasse et~al.(2019)Gasse, Ch{\'e}telat, Ferroni, Charlin, and
  Lodi}]{r_18}
Gasse, M.; Ch{\'e}telat, D.; Ferroni, N.; Charlin, L.; and Lodi, A. 2019.
\newblock Exact combinatorial optimization with graph convolutional neural
  networks.
\newblock In \emph{Advances in Neural Information Processing Systems}.

\bibitem[{Gilmore and Gomory(1961)}]{r_9}
Gilmore, P.~C.; and Gomory, R.~E. 1961.
\newblock A linear programming approach to the cutting-stock problem.
\newblock \emph{Operations Research}, 9(6): 849--859.

\bibitem[{Goffin and Vial(2000)}]{r_6}
Goffin, J.-L.; and Vial, J.-P. 2000.
\newblock Multiple cuts in the analytic center cutting plane method.
\newblock \emph{SIAM Journal on Optimization}, 11(1): 266--288.

\bibitem[{{Gurobi Optimization, LLC}(2023)}]{r_36}
{Gurobi Optimization, LLC}. 2023.
\newblock {Gurobi optimizer reference manual}.

\bibitem[{Kami{\'n}ski et~al.(2021)Kami{\'n}ski, Ludwiczak, Jasi{\'n}ski,
  Bukala, Madaj, Szczepaniak, and Dunin-Horkawicz}]{r_21}
Kami{\'n}ski, K.; Ludwiczak, J.; Jasi{\'n}ski, M.; Bukala, A.; Madaj, R.;
  Szczepaniak, K.; and Dunin-Horkawicz, S. 2021.
\newblock {Rossmann-toolbox: A deep learning-based protocol for the prediction
  and design of cofactor specificity in Rossmann fold proteins}.
\newblock \emph{Briefings in Bioinformatics}, 23(1): bbab371.

\bibitem[{L{\"{u}}bbecke and Desrosiers(2005)}]{r_1}
L{\"{u}}bbecke, M.~E.; and Desrosiers, J. 2005.
\newblock Selected topics in column generation.
\newblock \emph{Operations Research}, 53(6): 1007--1023.

\bibitem[{Malaguti and Toth(2010)}]{r_22}
Malaguti, E.; and Toth, P. 2010.
\newblock A survey on vertex coloring problems.
\newblock \emph{International Transactions in Operational Research}, 17(1):
  1--34.

\bibitem[{Mazyavkina et~al.(2021)Mazyavkina, Sviridov, Ivanov, and
  Burnaev}]{r_7}
Mazyavkina, N.; Sviridov, S.; Ivanov, S.; and Burnaev, E. 2021.
\newblock Reinforcement learning for combinatorial optimization: {A} survey.
\newblock \emph{Computers \& Operations Research}, 134: 105400.

\bibitem[{Mehrotra and Trick(1996)}]{r_10}
Mehrotra, A.; and Trick, M.~A. 1996.
\newblock A column generation approach for graph coloring.
\newblock \emph{INFORMS Journal on Computing}, 8(4): 344--354.

\bibitem[{Morabit, Desaulniers, and Lodi(2021)}]{r_16}
Morabit, M.; Desaulniers, G.; and Lodi, A. 2021.
\newblock Machine-learning-based column selection for column generation.
\newblock \emph{Transportation Science}, 55(4): 815--831.

\bibitem[{Moungla, L{\'{e}}tocart, and Nagih(2010)}]{r_3}
Moungla, N.~T.; L{\'{e}}tocart, L.; and Nagih, A. 2010.
\newblock Solutions diversification in a column generation algorithm.
\newblock \emph{Algorithmic Operations Research}, 5(2): 86--95.

\bibitem[{Pessoa et~al.(2018)Pessoa, Sadykov, Uchoa, and Vanderbeck}]{r_12}
Pessoa, A.; Sadykov, R.; Uchoa, E.; and Vanderbeck, F. 2018.
\newblock Automation and combination of linear-programming based stabilization
  techniques in column generation.
\newblock \emph{INFORMS Journal on Computing}, 30(2): 339--360.

\bibitem[{Schulman et~al.(2017)Schulman, Wolski, Dhariwal, Radford, and
  Klimov}]{r_8}
Schulman, J.; Wolski, F.; Dhariwal, P.; Radford, A.; and Klimov, O. 2017.
\newblock Proximal policy optimization algorithms.
\newblock \emph{arXiv preprint arXiv:1707.06347}.

\bibitem[{Shen et~al.(2022)Shen, Sun, Li, Eberhard, and Ernst}]{r_13}
Shen, Y.; Sun, Y.; Li, X.; Eberhard, A.; and Ernst, A. 2022.
\newblock Enhancing column generation by a machine-learning-based pricing
  heuristic for graph coloring.
\newblock In \emph{AAAI Conference on Artificial Intelligence}.

\bibitem[{Tahir et~al.(2021)Tahir, Quesnel, Desaulniers, Hallaoui, and
  Yaakoubi}]{r_15}
Tahir, A.; Quesnel, F.; Desaulniers, G.; Hallaoui, I.~E.; and Yaakoubi, Y.
  2021.
\newblock An improved integral column generation algorithm using machine
  learning for aircrew pairing.
\newblock \emph{Transportation Science}, 55(6): 1411--1429.

\bibitem[{Vanderbeck(1994)}]{r_4}
Vanderbeck, F. 1994.
\newblock \emph{Decomposition and column generation for integer programs}.
\newblock Ph.D. thesis, UCL-Universit{\'e} Catholique de Louvain.

\bibitem[{Veli{\v{c}}kovi{\'{c}} et~al.(2018)Veli{\v{c}}kovi{\'{c}}, Cucurull,
  Casanova, Romero, Li{\`{o}}, and Bengio}]{r_20}
Veli{\v{c}}kovi{\'{c}}, P.; Cucurull, G.; Casanova, A.; Romero, A.; Li{\`{o}},
  P.; and Bengio, Y. 2018.
\newblock {Graph attention networks}.
\newblock In \emph{International Conference on Learning Representations}.

\bibitem[{Xu et~al.(2019)Xu, Hu, Leskovec, and Jegelka}]{r_19}
Xu, K.; Hu, W.; Leskovec, J.; and Jegelka, S. 2019.
\newblock How powerful are graph neural networks?
\newblock In \emph{International Conference on Learning Representations}.

\bibitem[{Yang, Jiang, and Song(2023)}]{r_37}
Yang, W.; Jiang, P.; and Song, S. 2023.
\newblock High-speed Train Timetabling Based on Reinforcement Learning.
\newblock In \emph{IEEE Symposium Series on Computational Intelligence},
  1187--1193.

\bibitem[{Yuan, Jiang, and Song(2022)}]{r_14}
Yuan, H.; Jiang, P.; and Song, S. 2022.
\newblock The neural-prediction based acceleration algorithm of column
  generation for graph-based set covering problems.
\newblock In \emph{IEEE International Conference on Systems, Man, and
  Cybernetics}, 1115--1120.

\bibitem[{Zhang et~al.(2020)Zhang, Song, Cao, Zhang, Tan, and Xu}]{r_38}
Zhang, C.; Song, W.; Cao, Z.; Zhang, J.; Tan, P.~S.; and Xu, C. 2020.
\newblock Learning to Dispatch for Job Shop Scheduling via Deep Reinforcement
  Learning.
\newblock In \emph{Advances in Neural Information Processing Systems}.

\end{thebibliography}

\newpage

\section*{Appendix}

\subsection*{A. Hardware and Environment}

The training and inference phases of our approach are implemented on NVIDIA GeForce RTX 3090 24G GPU and based on PyTorch 1.12.1 and Deep Graph Library 0.9.1. The RMP and PP in CG are solved using Gurobi 9.5.2. All optimization procedures are run on AMD Ryzen ThreadRipper 3990X CPU @ 2.90GHz.

\subsection*{B. Features Used for the State Representation}

Here we provide details of the features used in the state representation of the MDP (see Figure \ref{fig_2}), which is embedded through the neural network for the multiple-column selection strategy.

\subsubsection*{Constraint node features}
\begin{itemize}
    \item  \textbf{Dual value}: Dual value (or shadow price) corresponding to each constraint, which is then used in the objective function of PP.
    \item  \textbf{Connectivity of constraint node}: The total number of column nodes each constraint node connects to, indicating the connectivity of each constraint node on the bipartite graph.
    \item  \textbf{Right-hand side value}: Right-hand side value of the corresponding constraint in RMP, which indicates the demand $d_i$ in CSP and is always equal to 1 in GCP.
    \item  \textbf{Slack in the current solution}: The constraint slack in the current RMP solution, which is equal to the left-hand side value according to the current solution minus the right-hand side value of the constraint.
\end{itemize}

\subsubsection*{Column node features}
\begin{itemize}
    \item  \textbf{Reduced cost}: A quantity associated with each variable, indicating how much the objective function coefficient would have to improve before it would be possible for a corresponding variable to assume a positive value in the optimal solution.
    \item  \textbf{Connectivity of column node}: The total number of constraint nodes each column node connects to, indicating the connectivity of each column node on the bipartite graph.
    \item  \textbf{Solution value}: The solution value of each decision variable corresponding to each column node after solving the RMP in the current iteration.
    \item  \textbf{Waste}: A feature recording the remaining length of a roll if the current pattern is cut from the roll (used in CSP).
    \item  \textbf{If the column left basis on the last iteration or not}: A binary feature with the value 1 if and only if the corresponding column has left basis at the last iteration.
    \item  \textbf{If the column entered basis on the last iteration or not}: A binary feature with the value 1 if and only if the corresponding column has entered basis at the last iteration.
    \item  \textbf{Number of iterations that a column node stays in the basis}: If a column node stays in the basis for a long time, it is likely that the variable corresponding to this column node is really good.
    \item  \textbf{Number of iterations that a column node stays out of the basis}: If a column node stays out of the basis for a long time, it is likely that it will never enter the basis and be used in the optimal solution in future iterations.
    \item  \textbf{Candidate column node or not}: A binary feature indicating whether the column node is a candidate or not. If the column node is a candidate column node generated at the current iteration by PP, then this binary feature takes 1 otherwise 0.
\end{itemize}

\subsubsection*{Global features}\,\par
The global features used for the CSP instances:
\begin{itemize}
    \item  \textbf{Roll length}: The length $L$ of the roll that has to be cut to satisfy the different demands.
    \item  \textbf{Total demands}: The sum of demands of the pieces ordered by all customers: $\sum_{i=1}^{m}d_i$.
    \item \textbf{Minimum length of demand pieces}: The minimum length of the pieces ordered by the customers, which is normalized by the roll length $L$: $\min_{i\in\left\{1,\ldots,m\right\}}{\ell_i/L}$.
    \item \textbf{Maximum length of demand pieces}: The maximum length of the pieces ordered by the customers, which is normalized by the roll length $L$: $\max_{i\in\left\{1,\ldots,m\right\}}{\ell_i/L}$.
\end{itemize}

The global features used for the GCP instances:
\begin{itemize}
    \item  \textbf{Number of nodes}: The number of nodes in the graph needed to be assigned colors such that every pair of adjacent nodes does not share the same color.
    \item  \textbf{Number of edges}: The total number of edges connecting the nodes in the graph. The value is normalized by the number of node pairs, indicating the possibility that an edge occurs between a node pair.
\end{itemize}

\begin{algorithm*}[t]
\caption{Proximal Policy Optimization}\label{alg_ppo}
\textbf{Input}: initial actor model parameters $\theta_0$, initial critic model parameters $\phi_0$.\\
\textbf{Hyper-parameter}: clipping threshold $\epsilon$.\\
\textbf{Output}: trained actor model parameters $\theta$, trained critic model parameters $\phi$.
\begin{algorithmic}
\FOR {$k = 0,1,2,...$}
\STATE Collect the set of trajectories $\mathcal{D}_{k}=\left\{\tau_{i}\right\}$ by running the actor policy $\pi_{\theta_{k}}$. \\
\STATE Compute the rewards-to-go $\hat{R}_{t}$. \\
\STATE Compute the advantage estimates $\hat{A}_{t}$ based on the current value function $V_{\phi_{k}}$ from the critic $\phi_{k}$. \\
\STATE Update the actor model parameters by maximizing the PPO-clip objective:\\
$$\theta_{k+1}=\arg \max _{\theta} \frac{1}{\left|\mathcal{D}_{k}\right|} \sum_{\tau \in \mathcal{D}_{k}}\left( \frac{1}{T_{\tau}} \sum_{t=0}^{T_{\tau}} \min \left(\frac{\pi_{\theta}\left(a_{t} \mid s_{t}\right)}{\pi_{\theta_{k}}\left(a_{t} \mid s_{t}\right)} \hat{A}_{t}, \quad \operatorname{clip}\left(\frac{\pi_{\theta}\left(a_{t} \mid s_{t}\right)}{\pi_{\theta_{k}}\left(a_{t} \mid s_{t}\right)}, 1-\epsilon, 1+\epsilon\right)\hat{A}_{t}\right) \right),$$\\
typically via stochastic gradient ascent with Adam.
\STATE Fit the value function by regression on the  mean-squared error (MSE):\\
$$\phi_{k+1}=\arg \min _{\phi} \frac{1}{\left|\mathcal{D}_{k}\right|} \sum_{\tau \in \mathcal{D}_{k}} \left( \frac{1}{T_{\tau}} \sum_{t=0}^{T_{\tau}}\left(V_{\phi}\left(s_{t}\right)-\hat{R}_{t}\right)^{2} \right)$$\\
typically via some gradient descent algorithm.\\
\ENDFOR
\end{algorithmic}
\end{algorithm*}

\subsection*{C. Data Generation}

\subsubsection*{CSP instance}\,

We follow \citet{r_24, r_23} in the generation for instances of roll length $L=50, 100,$ and $200$ as the configuration in Table \ref{table_5}.

According to the configuration, $n$ pieces of demand are generated independently, each of which gets a length uniformly distributed in $\left[w_{\text{min}}L, w_{\text{max}}L \right]$. All pieces of the same length are merged as one order.

The CG formulation of CSP is introduced in the previous section.

\begin{table}[h]
\renewcommand\arraystretch{1.3}
    \setlength{\tabcolsep}{4pt}
    \centering
    \caption{Instance configuration for CSP.}
    \fontsize{9pt}{9pt}\selectfont
    \label{table_5}
    \begin{tabular}{c|cccc}
        \toprule
        Instance & $L$ & $n$ & $w_{\text{min}}$ & $w_{\text{max}}$ \\
        \midrule
        Easy & $50$ & $\left\{50, 75, 100, 120\right\}$ & $\left\{0.1, 0.2\right\}$ & $\left\{0.7, 0.8\right\}$ \\
        Normal & $100$ & $\left\{575, 100, 120, 150\right\}$ & $\left\{0.1, 0.2\right\}$ & $\left\{0.7, 0.8\right\}$ \\
        Hard & $200$ & $\left\{125, 150\right\}$ & $\left\{0.1, 0.2\right\}$ & $\left\{0.7, 0.8\right\}$ \\
        \bottomrule
    \end{tabular}
\end{table}

\subsubsection*{GCP instance}\,

We follow \citet{r_10} in the generation for instances of the number of nodes $N=30, 40,$ and $50$, and each edge occurs independently with a probability $p$, where $p \in \left[0.4, 0.6\right]$.

The detailed CG formulation of GCP is as follows. Let $G(\mathcal{V}, \mathcal{E})$ denotes a graph, where $\mathcal{V}$ is the node set and $\mathcal{E}$ is the edge set, then the GCP can be formulated as a set partitioning problem:

\begin{align*}
    \text{ min } & \sum_{s \in \mathcal{S}} x_{s}\\
    \text { s.t. } & \sum_{s \in \mathcal{S}} \mathbb{I}\left(i \in s \right) x_{s} = 1, \; i \in \mathcal{V}, \\
    & x_{s} \in \left\{0, 1\right\}, \; s \in \mathcal{S},
\end{align*}
where $\mathcal{S}$ is the set of all the possible MISs in that graph. The binary variable $x_s$ indicates whether an MIS $s$ is used to cover a graph, and the indicator function $\mathbb{I}\left(i \in s \right)$ takes value $1$ if node $i$ is in MIS $s$ and value $0$ otherwise.

The RMP is a linear relaxation with a fraction of the MISs $\tilde{\mathcal{S}} \subset \mathcal{S}$:

\begin{align*}
    \text{ min } & \sum_{s \in \tilde{\mathcal{S}}} x_{s}\\
    \text { s.t. } & \sum_{s \in \tilde{\mathcal{S}}} \mathbb{I}\left(i \in s \right) x_{s} = 1, \; i \in \mathcal{V}, \\
    & 0 \le x_{s} \le 1, \; s \in \tilde{\mathcal{S},}
\end{align*}
while the PP is modeled as a maximum weight independent set problem (MWISP) to search for new MISs with the most negative reduced cost:

\begin{align*}
    \text{ max } & \sum_{i \in \mathcal{V}} \pi_{i}v_i \\
    \text { s.t. } & v_i + v_j \leq 1, \; \left(i,j\right) \in \mathcal{E}, \\
    & v_{i} \in \left\{0, 1\right\}, \; i \in \mathcal{V},
\end{align*}
where the binary variable $v_i$ denotes whether node $i$ is a part of the solution, i.e., an MIS according to the constraints. The dual value $\pi_{i}$ is the weight of the node $i$ in the MWISP.

\subsection*{D. PPO training}

We use PPO \cite{r_8} as the training algorithm for the multiple-column selection strategy. PPO is a deep reinforcement learning algorithm based on the actor-critic architecture. Taking a state $s$ as input, the critic estimates the value function $V(s)$, and the actor gives a probability distribution (i.e., the policy) $\pi = (\pi\left(a_1 \mid s\right),$ $\pi\left(a_2\mid s\right), \dots, \pi\left(a_{\left | \mathcal{A} \right |}\mid s\right))$ over the action space. An action is sampled from the probability distribution, and the corresponding $k$ columns are selected and added to the RMP to start a new iteration.

The idea behind PPO is to improve the training stability of the actor policy by limiting the change to the policy at each training epoch. PPO constrains the actor policy update with a new objective function called the clipped surrogate objective function between $\left[1-\epsilon, 1+\epsilon\right]$, which clips the estimated advantage function if the new policy is far away from the old policy, to ensure the policy change in a small range. The detailed training process of PPO is presented in Algorithm \ref{alg_ppo}. 

\subsection*{E. Diverse-M Strategy}

We provide the implementation details of the Diverse-M strategy in our experiment. The idea behind Diverse-M is to select more diversiﬁed columns at each iteration. Diversified columns should not cover the same tasks as much as possible. Two columns that do not cover the same tasks are said to be \emph{disjoint}.

First, the candidate columns are sorted in increasing order of their reduced cost. Second, following this order, the candidate columns are distributed in blocks (numbered $1, 2, 3, ...$) of disjoint columns using the following rule: a column is put in block numbered $b$ if this column is not disjoint from at least one column in each block $1$ to $b-1$, and disjoint from all columns already in block $b$. Consequently, the first column is put in block $1$, and then the second is put in block $1$ if it is disjoint from the first column or in block $2$ otherwise, and so on. Then, the candidate columns in block $1$ are selected, and then block $2$, until the selected columns reach the limit.

\end{document}